\documentclass[fleqn]{amsart}
\usepackage{amssymb, graphics, amsmath, epsfig, changebar}

\theoremstyle{plain}
\newtheorem{theorem}{Theorem}
\newtheorem{corollary}{Corollary}

\newtheorem{lemma}{Lemma}
\newtheorem{proposition}{Proposition}

\theoremstyle{definition}
\newtheorem{definition}{\rm{Definition}}

\theoremstyle{remark}

\numberwithin{equation}{section}
\newcommand{\N}{\mathbb{N}}
\newcommand{\Z}{\mathbb{Z}}
\newcommand{\R}{\mathbb{R}}
\newcommand{\Sum}{\displaystyle\sum}
\newcommand{\Bigoplus}{\displaystyle\bigoplus}
\newcommand{\KB}{\mathcal{S}_{2,\infty}}
\newcommand{\RR}{\R P^3\sharp\R P^3}
\newcommand{\Om}{\Omega}

\begin{document}                   

\title[KBSM of $\RR$]{Kauffman 
bracket skein module of the connected sum of two projective spaces}
\author{Maciej Mroczkowski}

\begin{abstract}
Diagrams and Reidemeister moves for links in a twisted $S^1$-bundle over 
an  unorientable surface are introduced. Using these diagrams, we compute 
the Kauffman Bracket Skein Module (KBSM) of $\RR$. In particular, 
we show that it has torsion. We also present a new computation of the 
KBSM of $S^1\times S^2$ and the lens spaces $L(p,1)$.
\end{abstract}
\maketitle

\let\thefootnote\relax 

\section{Introduction}
Skein modules, which are invariants of 3-manifolds as well as of links in 
these manifolds, were introduced by J. Przytycki~\cite{P1} and V. 
Turaev~\cite{T}. In this paper we compute the Kauffman bracket skein 
module (KBSM) of $\RR$ (Theorem \ref{main_theorem}). This is the first 
such computation for a closed non-prime manifold. 
Also, it is the second example of a fully computed KBSM with torsion  
for a closed manifold, after the case of $S^1\times S^2$~\cite{HP}. Unlike 
the KBSM of $S^1\times S^2$, the KBSM of $\RR$ does not split as a sum of 
cyclic modules (Proposition \ref{non_split}).

\begin{figure}[h]
\scalebox{0.8}{\includegraphics{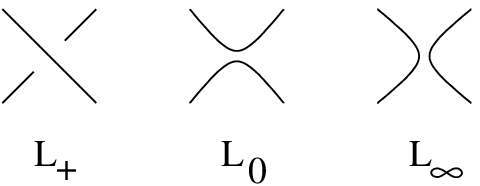}}
\caption{}
\label{smooth}
\end{figure}

First, we recall the definition of this module.
Throughout this paper $R$ will be the ring of Laurent polynomials in $A$, 
$R=\mathbb Z[A,A^{-1}]$. Let $M$ be an orientable 3-manifold.
The Kauffman bracket skein module of $M$, or $\KB(M)$, is the $R$-module
generated by isotopy classes of unoriented framed links in $M$ modulo
local relations:\\ \\
$(K1):\; L_{+}=A L_0 + A^{-1} L_{\infty}$\\
$(K2):\; L\sqcup T=(-A^2-A^{-2}) L$\\ \\
where $T$ is the trivial framed knot and the triple $L_{+}$, $L_0$ and 
$L_\infty$ is presented in Figure~\ref{smooth}.

It can be seen easily that $L^{(1)}=-A^{3}L$ in 
$\mathcal{S}_{2,\infty}(M)$ where $L^{(1)}$ is obtained from $L$ by adding 
a positive twist. It is called the {\it framing relation}.

For example $\KB(S^3)$ is free cyclic, $\KB(F \times I)$ ($F$ orientable 
surface) is free, generated by isotopy classes of simple closed curves on 
$F$ without trivial components~\cite{P1}, $\KB(L(p,q))$ is free with 
$\lfloor p/2\rfloor +1$ generators~\cite{HP2}, whereas $\KB(S^1\times 
S^2)$ has torsion~\cite{HP}.

In the next section we introduce diagrams and Reidemeister moves for links in 
$N\hat{\times}S^1$, where $N$ is an unorientable surface. In section 3
we recall some results about $\KB(S^1\times B^2)$, $B^2$ a disk, 
from~\cite{MD}. In section 4 we compute $\KB(\RR)$, we show that it 
contains torsion elements and does not split as a sum of cylic 
$R$-modules. Finally, in sections 5 and 6 we present a new way to compute 
$\KB(S^1\times S^2)$ and $\KB(L(p,1))$.

\section{Diagrams of links in $N\hat{\times} S^1$ and Reidemeister moves}

In \cite{MD} the notion of diagrams of links in $F\times S^1$, where $F$ 
is an orientable surface, was introduced. Also, the Redeimeister moves for 
such diagrams were found. We recall these notions here, as they are 
essential for our construction. 

To obtain a diagram of a link $L$ in $F\times S^1$, we cut $F\times S^1$ 
along $F\times \{1\}$, $1\in S^1$, thus obtaining $F\times [0,1]$ in 
which $L$ becomes $L'$, a collection of arcs. Then, $L'$ is projected 
onto $F\times \{0\}\approx F$ via a vertical projection yielding a set of 
closed curves in $F$ on which we keep some extra information: for each 
double point the usual information of over- and undercrossing depending 
on the relative height in $[0,1]$ of the projected points; for points 
coming from endpoints of arcs in $L'$ (which appeared after cutting 
$F\times S^1$) a dot with an arrow on it indicating the direction of 
increasing height in $[0,1]$ just before and just after the cut. In other 
words, travelling on $L$ in the direction of the arrow, and crossing it, 
one crosses the "roof" $F\times \{1\}$ and one emerges from the "floor" 
$F\times \{0\}$ in the cylinder $F\times [0,1]$.

\begin{figure}[h]
\scalebox{0.8}{\includegraphics{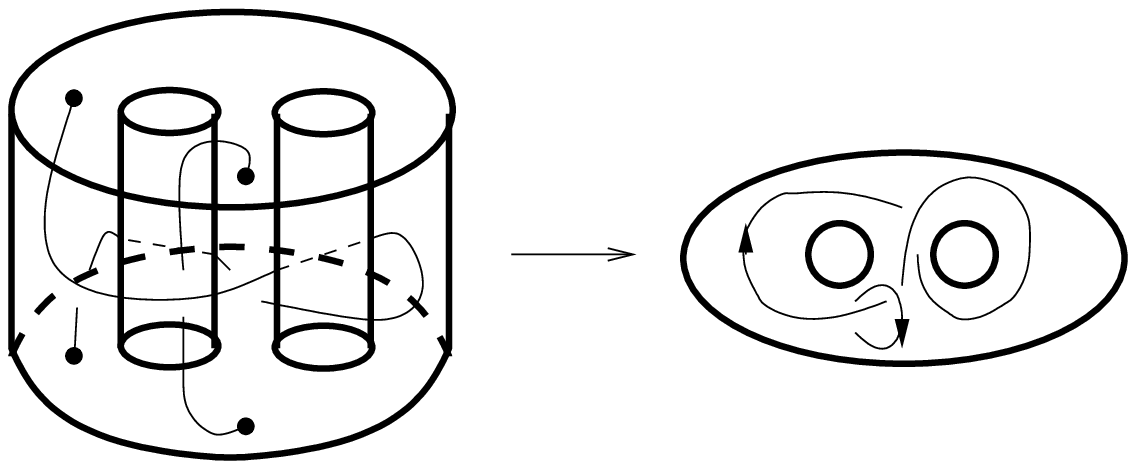}}
\caption{}
\label{diag} 
\end{figure}

The construction of a diagram in the case where $F$ is a disk with two 
holes is pictured in Figure~\ref{diag}.

\begin{figure}[h]
\scalebox{0.8}{\includegraphics{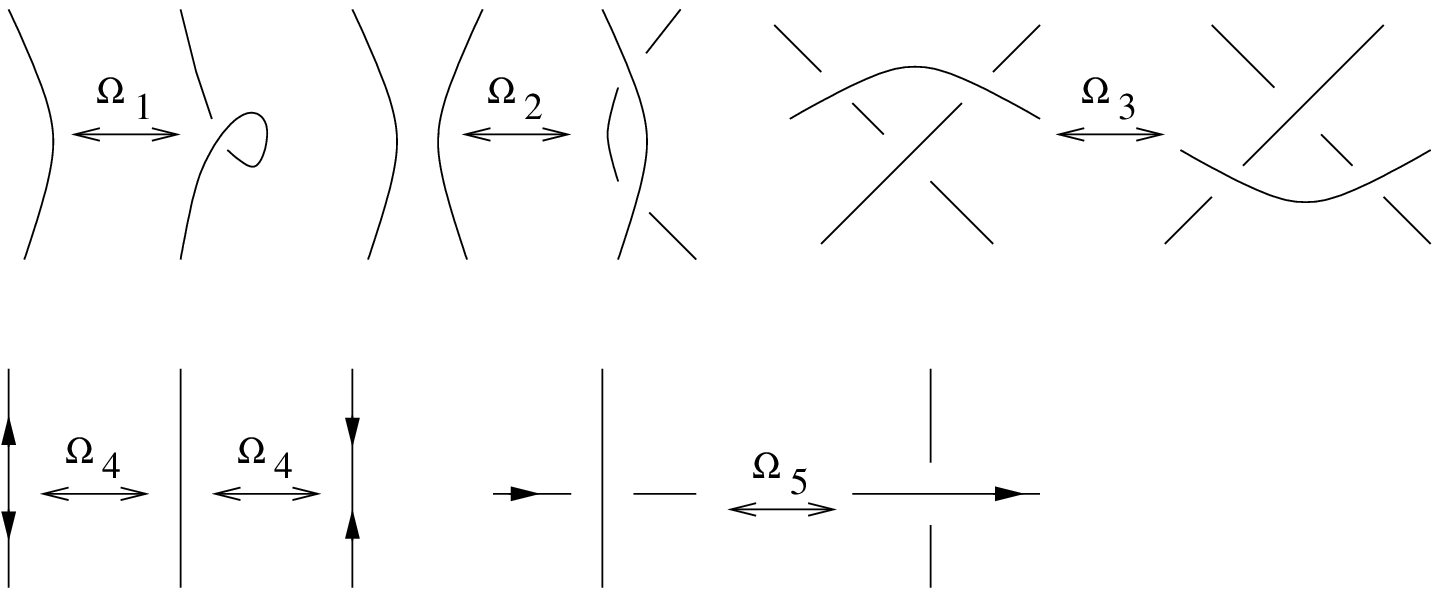}}
\caption{}
\label{reid} 
\end{figure}

\begin{figure}[h]
\scalebox{0.8}{\includegraphics{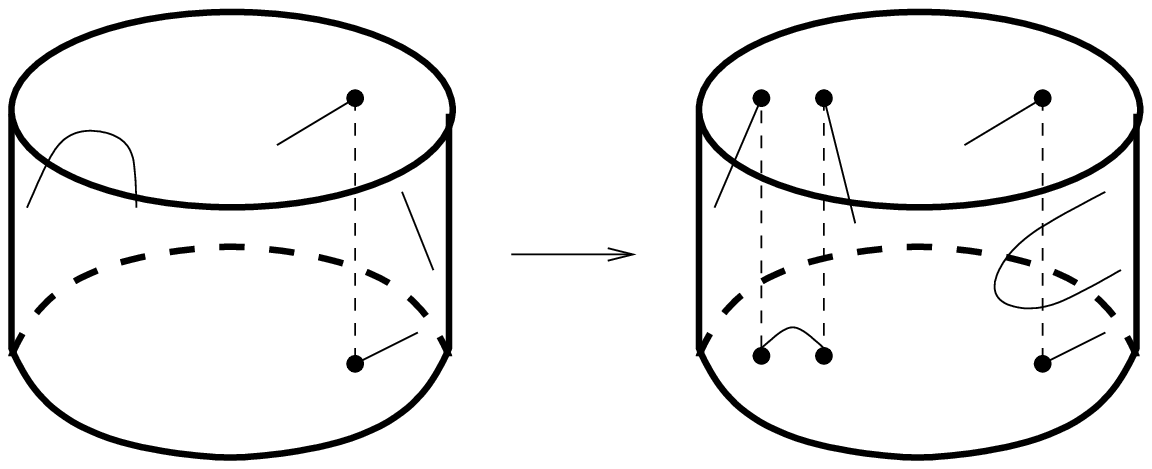}}
\caption{}
\label{interpr45} 
\end{figure}

Two links in $F\times S^1$ are isotopic if one can get from any diagram of 
one to any diagram of the other through a series of five Reidemeister 
moves, pictured in Figure~\ref{reid}. 
The interpretation of $\Om_4$ and $\Om_5$ is pictured 
in Figure~\ref{interpr45}.

We introduce now diagrams and Reidemeister moves for $N\hat{\times}S^1$, 
where $N$ is an unorientable surface. Denote an unorientable surface of
genus $k$ and $n-k$ boundary components by $N_{n,k}$. 
Such $N_{n,k}$ is constructed from a 2-sphere with $n$ holes, denoted 
$S^2_n$, by glueing M\"obius bands to $k$ ($k\le n$) of these holes, or
equivalently, by glueing $k$ boundary components ($S^1$-s) of $S^2_n$,
each component being glued to itself via the antipodal map. Denote these
$k$ components  by $C$. Denote the image after the glueing of $C$ in 
$N_{n,k}$ by $C'$.

\begin{figure}[h]
\scalebox{0.6}{\includegraphics{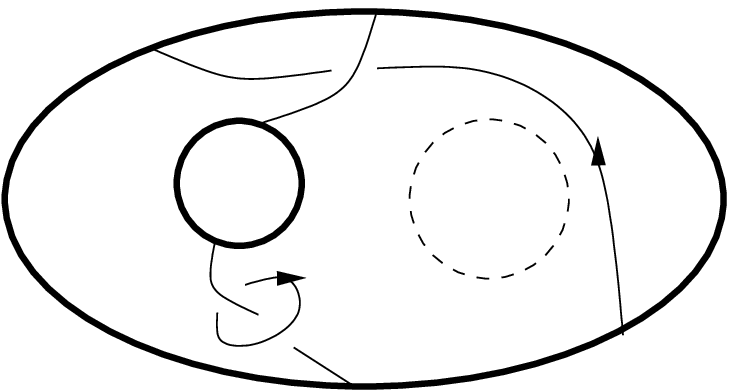}}
\caption{}
\label{crosscaps} 
\end{figure}

Let $L$ be a link in $N\hat{\times} S^1$ where $N=N_{n,k}$ for some 
$k\le n$. Cutting $N$ along $C'$ gives $S^2_n$ and cutting $N\hat{\times} 
S^1$ along $C'\hat{\times} S^1$ gives $S^2_n\times S^1$, in which $L$ 
becomes $L'$, a collection of arcs with endpoints in $\partial 
S^2_n\times S^1$. For such $L'$ in $S^2_n\times S^1$ a diagram is 
constructed in the same way as it was done for links in $F\times S^1$, 
with orientable $F$.
A diagram with double points (equipped with the information of 
over- and undercrossing) and dots with arrows is thus obtained. The 
difference is that now it may consist of arcs (and not only of closed 
curves), where the endpoints of arcs lie in $C$ in antipodal pairs. 
An example of a diagram of a link in $N_{3,2}\hat{\times}S^1$ is shown
in Figure~\ref{crosscaps}.

\begin{figure}[h]
\scalebox{0.8}{\includegraphics{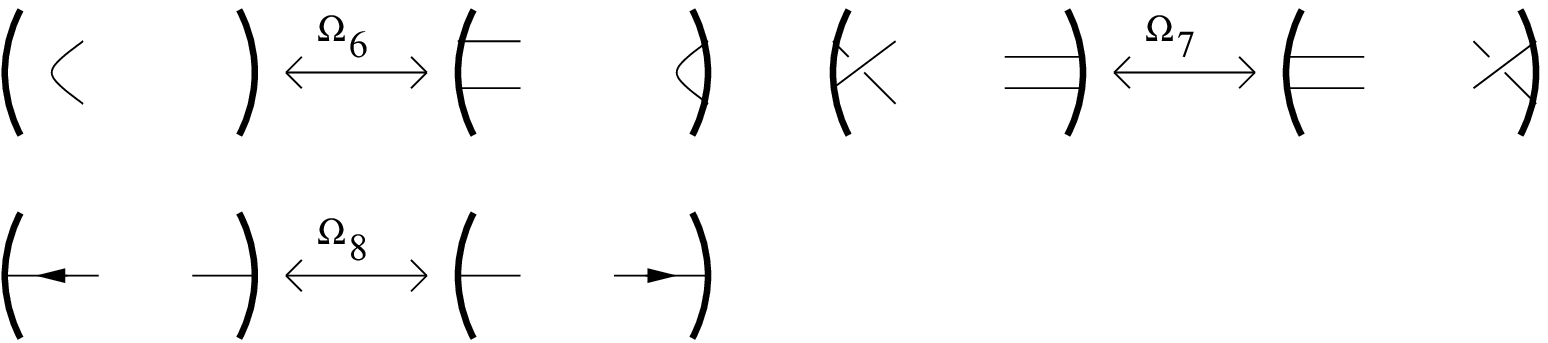}}
\caption{}
\label{reid678} 
\end{figure}

As before, there are five Reidemeister moves (three classical and two new ones 
involving the arrows). Moreover, there are extra moves obtained by considering 
the resolution of generic singularities for the diagrams.
Three new generic singularities are possible: an arc can be tangent to 
$C$, a double point can lie in $C$ or a dot with an arrow can lie in $C$.
Resolving these singularities gives respectively $\Om_6$, $\Om_7$ and $\Om_8$ 
moves pictured in Figure~\ref{reid678}.

In $\Om_7$ the upper branch on one side becomes the lower branch after 
going through $C$ as the $S^1$-bundle is twisted on $C$. For the same 
reason in $\Om_8$ the arrow switches orientation after going through 
$C$.

The {\it regular} Reidemeister moves are all the Reidemeister moves except $\Om_1$.

In the rest of this paper we will consider links and their diagrams 
in $N=N_{1,1}=\R P^2$, so that $N\hat{\times}S^1$ is $\RR$.

\section{The module $\KB(S^1\times B^2)$}

In this section we recall some results from \cite{MD} concerning 
$\KB(S^1\times B^2)$, $B^2$ a disk, as they are used later. In the 
preceding section we 
described diagrams and Reidemeister moves for $F\times S^1$, where $F$ is 
any orientable surface, so in particular for $B^2 \times S^1$.

Let $\epsfig{file=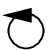}$ be denoted by $x$. 
Thus, the diagram $x$ represents a knot that runs parallel to the $S^1$ 
core of $S^1\times B^2$.
Applying the framing relation (i.e. $\Om_1$), $\Omega_5$, and another 
framing relation one gets the following easy Lemma (\cite{MD} Lemma 3.3):

\begin{lemma}\label{lemma_pm}
In $\KB(S^1\times B^2)$, 
$\epsfig{file=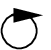}=A^{-6}\epsfig{file=x.eps}$ 
\end{lemma}

\begin{figure}[h]
\scalebox{0.8}{\includegraphics{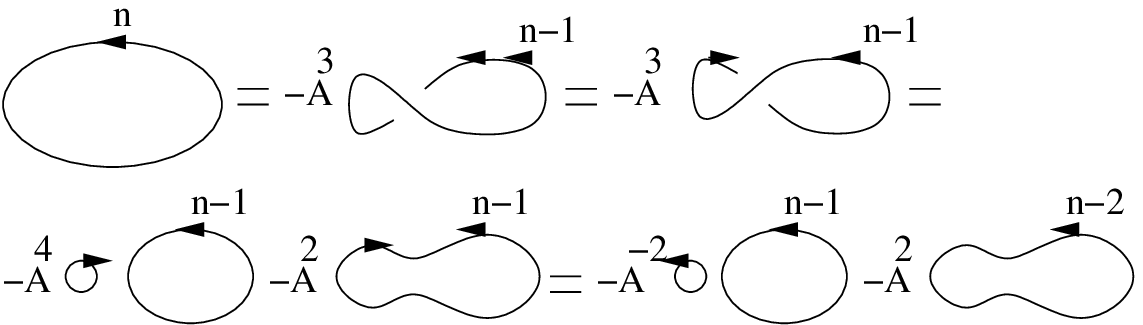}}
\caption{}
\label{D_n}
\end{figure}

Let $D_n$ be the diagram with no crossings and one component with $n$ 
arrows on it, $n\in\mathbb Z$ (if $n>0$ the arrows are counterclockwise 
and if $n<0$ the arrows are clockwise). Figure~\ref{D_n} shows how to 
express it, as an element of $\KB(S^1\times B^2)$, with $D_{n-1}$ and 
$D_{n-2}$. The last equality comes from Lemma \ref{lemma_pm}.

So, in the skein module we have the relation:
\begin{equation*}
D_n=-A^{-2}x D_{n-1}-A^2 D_{n-2}\label{eq1}
\end{equation*}

This suggests the following definition: 

\begin{definition}\label{P_n}
Let $P_n$, $n\in\mathbb Z$, be polynomials in $x$ with 
coefficients in the ring $R$ defined inductively by:\\ 
$P_0=-A^2-A^{-2}$\\ 
$P_1=x$\\ 
$P_n=-A^{-2}x P_{n-1}-A^2 P_{n-2}$\\ where 
the last relation is also used to define $P_n$ for all negative $n$.\\ 
\end{definition}

\begin{proposition}\label{DS1}(\cite{MD} Proposition 3.7)
$\KB(S^1\times B^2)$ is a free $R$-module, generated by
$\{x^n|n\in\mathbb N\}$, where $x^n$ stands for $\epsfig{file=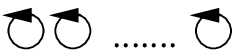}$
($n$ copies of $x$), $x^0=\emptyset$.
\end{proposition}

The isomorphism in this Proposition was given by a refined Kauffman bracket,
$<>_r$, which to a diagram $D$ associates its unique linear expression in 
$x^n$-s. This bracket respects the Kauffman relation $(K1)$ and $(K2)$ and
it was proven that it is invariant under regular Reidemeister moves.
We will use $<>_r$ in the next section. 

The bracket $<>_r$ was defined in such a way that it satisfies 
$<D_n>_r=P_n$, $n\in\Z$. Therefore, in what follows we will identify the 
polynomials $P_n$ with diagrams $D_n$, keeping in mind that they are 
polynomials in $x$ in $\KB(S^1\times B^2)$.

\section{The module $\KB(\RR)$}

\begin{figure}[h]
\scalebox{0.6}{\includegraphics{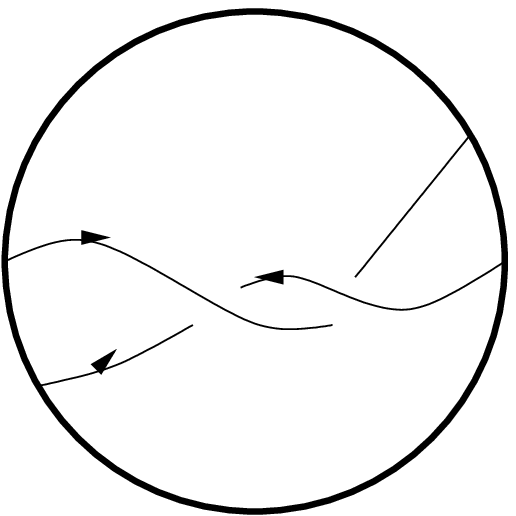}}
\caption{}
\label{ex_diag}
\end{figure}

In this section, we compute the KBSM of $\RR$, using diagrams of links in 
a twisted $S^1$-bundle over $\R P^2$, which is the same as $\RR$.
Thus, a diagram will consist of a 
family of arcs in a disc $B^2$, with endpoints lying on $\partial B^2$ 
in antipodal couples, information of under- and overcrossing for each 
double point, and with some arrows on the arcs outside the 
endpoints and outside the double points.
An example of such a diagram is pictured in Figure~\ref{ex_diag}.

\begin{figure}[h]
\scalebox{0.8}{\includegraphics{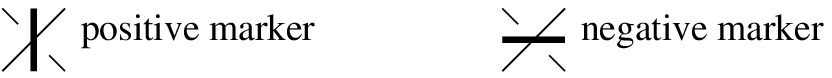}}
\caption{}
\label{markers}
\end{figure}

Let $D$ be a diagram of a link in $\RR$. In the same way as in the 
construction of the classical Kauffman bracket, each crossing of $D$ can 
be equipped with a positive or negative marker as is pictured in 
Figure~\ref{markers}.

A state $s$ of $D$ is the choice of a marker for each crossing of $D$. 
Let $D(s)$ be the diagram obtained from 
$D$ by smoothing all crossings in $D$ according to the markers 
determined by $s$. 

$D(s)$, viewed as a family of curves in $\R P^2$ (by identifying 
antipodal boundary points of the disk of the diagram to get $\R P^2$), 
consists of ovals and, at most, one projective line. The ovals may be 
nested (an oval is {\it nested} in another one, if it lies in the disk 
that this second oval bounds).
To each oval is associated its {\it arrow number}: it is the 
algebraic number of arrows on this oval with an arrow being {\it 
positive} (resp. {\it negative}), if it gives a counterclockwise (resp. 
clockwise) orientation to the disk the oval bounds.
To the projective line is associated its {\it parity} $0$ or $1$ which 
is the parity of the number of arrows on it.
An oval is {\it trivial} if its arrow number is equal to $0$ and if all 
ovals nested in it (there may be no such ovals) have arrow numbers 
equal to $0$.

Let $D'(s)$ be the diagram consisting of all non trivial ovals of $D(s)$ 
and the projective line (if there is one in $D(s)$), which are arranged 
in such a way that the projective line consists of one arc  
(with $0$ or $1$ arrow on it, equal to the parity of the projective 
line) and the ovals consist of simple closed curves in the disk of the 
diagram, nested in the same manner as in $D(s)$, each oval having $n$ 
arrows on it, where $n$ is the arrow number of this oval, arranged in a 
counterclockwise or clockwise way depending on the sign of $n$. 
Furthermore, all the ovals in $D'(s)$ lie on one side of 
the projective line and, if the parity of this line is $1$, they lie 
above the projective line on which the arrow points to the right.

\begin{figure}[h]
\scalebox{0.6}{\includegraphics{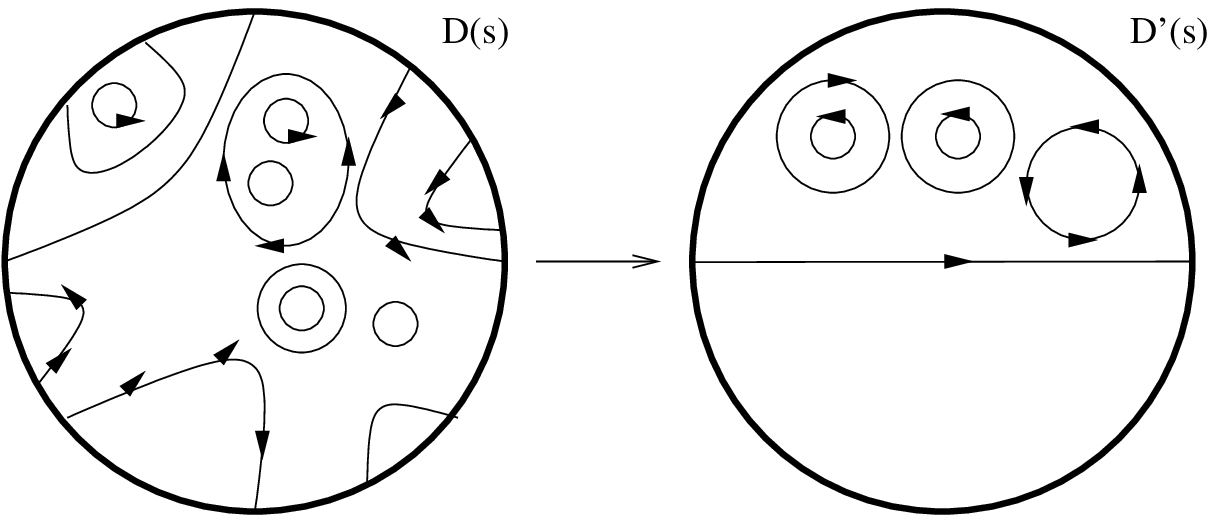}}
\caption{}
\label{DprimeS}
\end{figure}

An example of getting $D'(s)$ from $D(s)$ is pictured in 
Figure~\ref{DprimeS}.

Let $|s|$ be the number of trivial ovals of $D(s)$.
Let $p(s)$ (resp. $n(s)$) be the number of crossings with positive
(resp. negative) markers in $s$.

\begin{definition}\label{bracket}
The {\it Kauffman bracket} of $D$ is given by the following sum taken over 
all states of $D$:
\begin{equation*}<D>=\Sum_s A^{p(s)-n(s)}(-A^2-A^{-2})^{|s|}D'(s)
\end{equation*}
\end{definition}

\begin{lemma}
The Kauffman bracket $<>$ is preserved under all regular Reidemeister 
moves with the exception of $\Om_5$.
\begin{proof}[\rm{Proof}]
The proof of the invariance of $<>$ under $\Om_2$ and $\Om_3$-moves is analogous
to the classical case (Kauffman bracket for classical diagrams~\cite{K}).
For the moves $\Om_4$, $\Om_6$ and $\Om_8$ it suffices to prove the invariance of
$<>$ for diagrams without crossings (i.e. for each $D(s)$ separately). But then
this invariance follows from the easy fact that these moves do not change the nesting
of the ovals, the arrow numbers of the ovals and the partity of the arrows on the
projective line. Finally, the invariance of $<>$ under $\Om_7$ follows from the
invariance under $\Om_6$, which has to be applied twice as is pictured in 
Figure~\ref{inv_7}.

\begin{figure}[h]
\scalebox{0.8}{\includegraphics{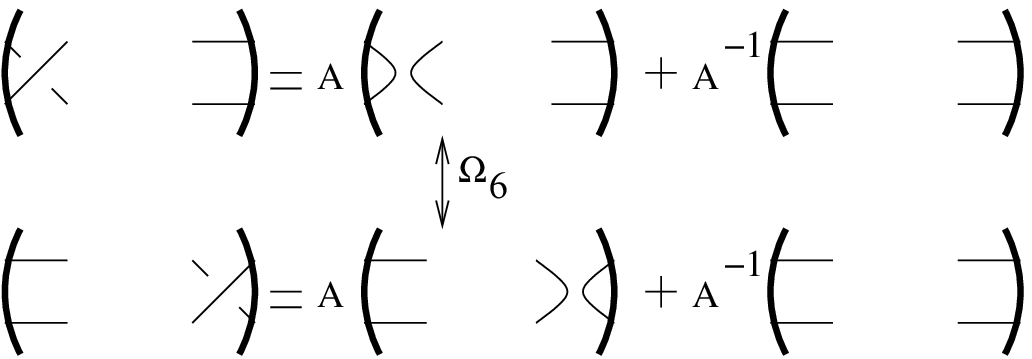}}
\caption{}
\label{inv_7}
\end{figure}

\end{proof}
\end{lemma}

We will refine now the bracket $<>$. Let $D$ be a diagram, $s$ a Kauffman
state of $D$. Then $D'(s)$ consists of ovals and, possibly, a projective 
line, as on the right of Figure~\ref{DprimeS}. All ovals lie in some disk 
$B$, which lies itself in the interior of the disk of the diagram. 
Therefore, they may be viewed as a diagram of a link in $S^1\times B$, 
for which a refined bracket $<>_r$ was defined in \cite{MD} using the 
isomorphism mentioned after Proposition \ref{DS1}. We use this bracket to 
define $<D'(s)>_r$  (outside of $B$ it leaves the projective 
line unchanged, if there is one present).

\begin{definition}\label{refined_bracket}
The {\it refined Kauffman bracket} of $D$ is given by the following sum 
taken over all states of $D$:
\begin{equation*}<D>_r=\Sum_s A^{p(s)-n(s)}(-A^2-A^{-2})^{|s|}<D'(s)>_r
\end{equation*}
\end{definition}

\begin{figure}[h]
\scalebox{0.6}{\includegraphics{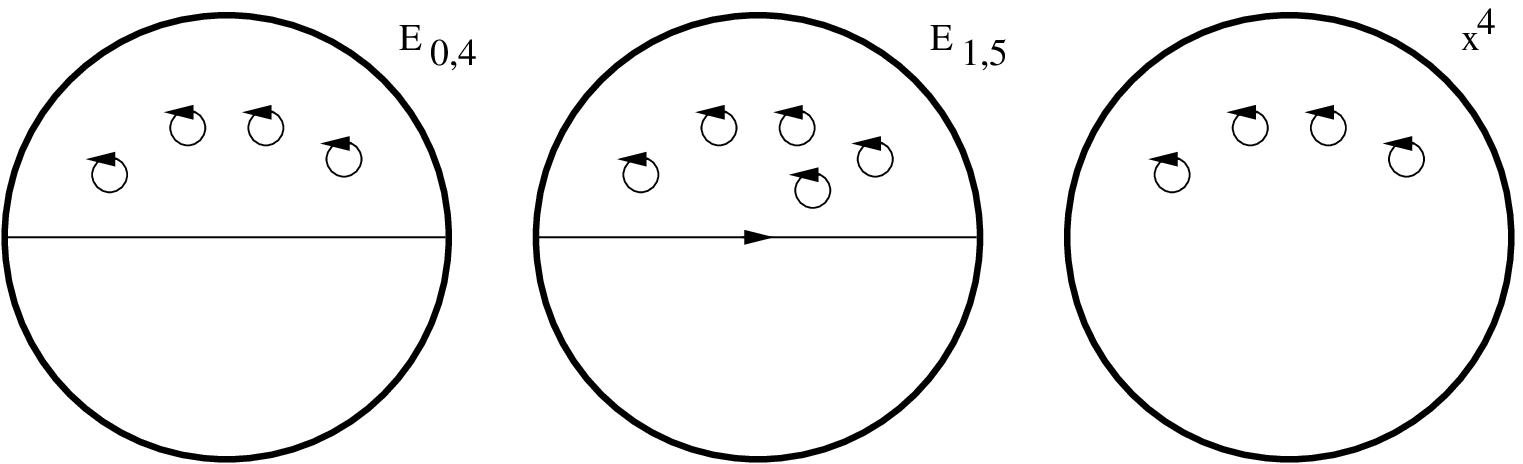}}
\caption{}
\label{basis_kbr}
\end{figure}

Thus, for a diagram $D$, $<D>_r$ is a linear expression in diagrams 
pictured in 
Figure~\ref{basis_kbr}. Denote the first two by $E_{m,n}$, where $m\in 
\Z_2$ is the
parity of the arrows on the projective line and $n$ is the number of $x$-s.

\begin{figure}[h]
\scalebox{0.8}{\includegraphics{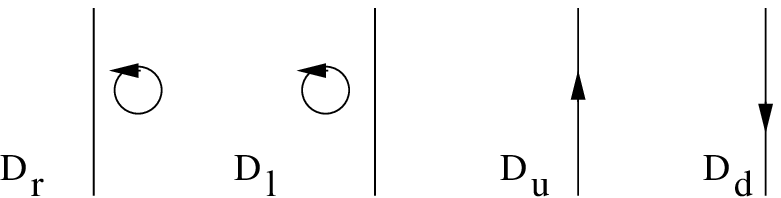}}
\caption{}
\label{xinout}
\end{figure}

Let $D_r$, $D_l$, $D_u$ and $D_d$ be four diagrams without crossings which differ only 
locally as pictured in Figure~\ref{xinout}. We assume also that the 
vertical strand
present in these diagrams belongs to an oval (i.e. not to a projective line).

Then we have:

\begin{lemma}\label{lemma_xinout}
The refined Kauffman bracket satisfies: 

(1) $<D_u>_r=-A^{-2}<D_r>_r-A^2<D_d>_r$

(2) $<D_u>_r=-A^{-4}<D_l>_r-A^{-2}<D_d>_r$

\begin{proof}[\rm{Proof}]
Note that by definition of $<>_r$, $<D_u>_r=<<D_u>>_r$ and the same is true for
the other three diagrams. Thus, we may replace each of the four diagrams with its
bracket (this reduces simply the number of endpoints of arcs on the 
boundary of the disk of the diagrams to $2$, if there is a projective 
line, or $0$ otherwise). 
The Lemma is now true because it is true for $<>_r$ defined on diagrams 
in $S^1\times B^2$, as was proven in \cite{MD} (Lemma 3.6), and we use 
here the same $<>_r$.
\end{proof}
\end{lemma}

The following Lemma will be used to reduce all the cases of $\Om_5$ moves to some
standard ones:

\begin{lemma}\label{lemma_standard5}
Suppose that $D_1$ is a diagram with one crossing and $D_2$ is obtained from $D_1$
with a single $\Om_5$ move. Then there is a series of $\Om_6$, $\Om_7$ and $\Om_8$
moves between $D_1$ and $D'_1$ and a similar series between $D_2$ and 
$D'_2$, where
$D'_2$ is obtained from $D'_1$ with a single $\Om_5$ move and $D'_1$ is 
one of the four types presented in Figure~\ref{standard5}. It is to be 
understood that 
in this Figure there are no other endpoints of arcs lying on the 
boundary of the disk, except, possibly, two such endpoints for types I 
and II corresponding to a projective line. 
There may be any configuration of ovals and arrows, including extra arrows 
on the arcs that are pictured.

\begin{figure}[h]
\scalebox{0.6}{\includegraphics{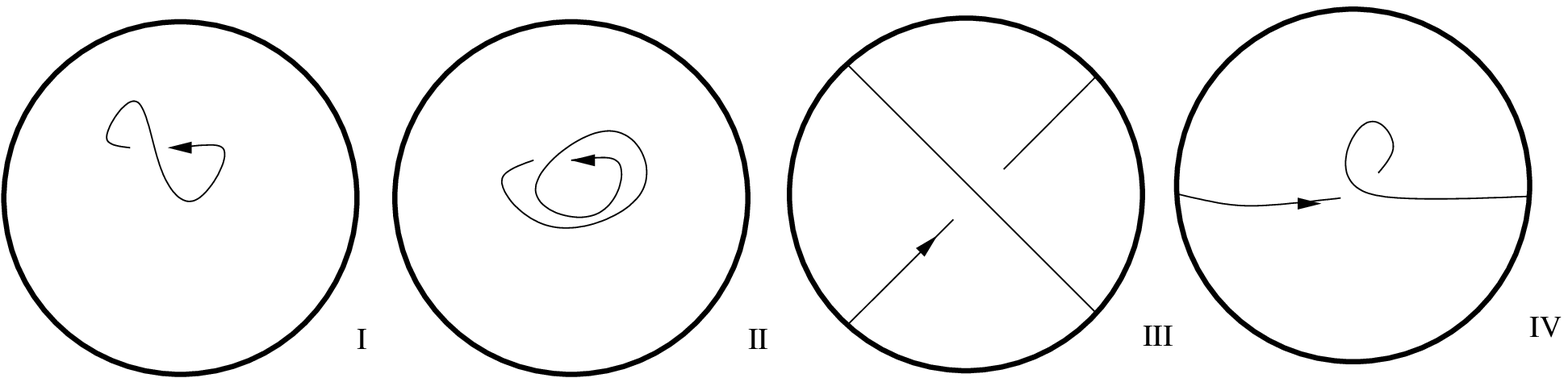}}
\caption{}
\label{standard5}
\end{figure}

\begin{proof}[\rm{Proof}]
The idea is to eliminate all arcs in $D_1$ and $D_2$ which do not 
contain the crossing and that have endpoints that are not antipodal.
Denote such an arc by $a$. Then $a$ divides the disk of diagrams $D_1$ 
and $D_2$ into two regions $R$ and $R'$ of which exactly one, say $R$, 
does not contain antipodal points on the boundary of the disk.
 
Now $a$ can be pushed through $R$ and across the boundary of the disk 
with the required 
moves after, possibly, having to push any other arcs and ovals that lie 
in $R$. It may happen that the unique crossing lies in $R$ in which 
case it is pushed across the boundary (both in $D_1$ and in $D_2$) 
giving again diagrams that differ by a single $\Om_5$ move.
This process reduces the number of endpoints of arcs on the boundary of 
the disk by at least two. Repeating it, one eliminates all such arcs. 

Now, if an endpoint belongs to an arc that does not contain the crossing, 
the other endpoint is antipodal to it and $D'_1$ is of type I or II.
Otherwise all endpoints belong to arcs that contain the crossing: if there 
are four such endpoints $D'_1$ is of type III, if there are two such 
endpoints, it is of type IV.
Finally if there are no endpoints of arcs, $D'_1$ is of type I or II.
\end{proof}
\end{lemma}

We will say that an $\Om_5$ move is of type I, II, III or IV, if it is 
a move of the type decribed between diagrams $D'_1$ and $D'_2$ in the 
hypothesis of the preceding Lemma. From the remarks after Proposition 
\ref{DS1}, $<>_r$ is invariant under moves of type I and II for diagrams 
of links in $S^1\times B^2$. As we use the same definition of $<>_r$ here, 
we also have:

\begin{lemma}\label{typesI_II}
The refined Kauffman bracket, $<>_r$, is invariant under $\Om_5$ moves of 
type I and II.
\end{lemma}

As $H_1(\RR;\Z_2)=\Z_2\oplus \Z_2$ and Kauffman relations $(K1)$ 
and $(K2)$ respect the $\Z_2$-homology, $\KB(\RR)$ splits into four 
submodules corresponding to each homology class: 
$(0,0)$, $(1,0)$, $(0,1)$ and $(1,1)$. We may therefore write:
\begin{equation*} \KB(\RR)=\KB^0(\RR)\oplus \KB^1(\RR) \end{equation*}
where $\KB^0(\RR)$ is generated by elements of classes $(0,0)$ and $(1,1)$
and $\KB^1(\RR)$ is generated by elements of classes $(1,0)$ and $(0,1)$.

There are two invariants for diagrams that do not change under Reidemeister moves
and Kauffman relations: the parity of the number of arrows and the parity 
of half the number of
the endpoints of arcs lying on the boundary of the disk.

If this second number is odd for $D$ (which is equivalent to the presence of a projective
line in $D(s)$ for all Kauffman states $s$), then the corresponding link represents
$(0,1)$ or $(1,0)$ in $H_1(\RR;\Z_2)$, the class depending on the parity of the arrows.
For example, the projective line can be viewed as lying in the first copy of $\R P^3$ in
$\RR$ and the projecive line with an arrow on it can be viewed as lying in the second
copy.

If half the number of the endpoints of arcs lying on the boundary of the 
disk is even (which is equivalent to the absence of the projective line in 
$D(s)$ for all
Kauffman states $s$), then the corresponding link represents $(0,0)$ or 
$(1,1)$ in $H_1(\RR;\Z_2)$, the class depending on the parity of the 
arrows. For example, the trivial knot realizes $(0,0)$, whereas a 
vertical $S^1$ in $\R P^2\hat{\times}S^1$ realizes $(1,1)$.

\subsection{The submodule $\KB^1(\RR)$}

\begin{figure}[h]
\scalebox{0.8}{\includegraphics{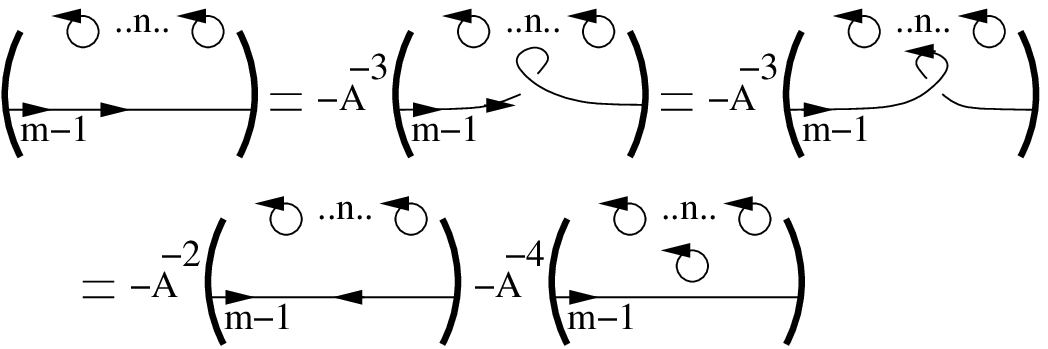}}
\caption{}
\label{red_x_5}
\end{figure}

A relation in $\KB^1(\RR)$ is pictured in Figure~\ref{red_x_5}.

We can express this relation using the notations introduced before 
(see Figure~\ref{basis_kbr}):
\begin{equation*} E_{m,n}=-A^{-2}E_{m-2,n}-A^{-4}E_{m-1,n+1}
\end{equation*}
Using that $m\in\Z_2$ and rearanging the terms gives:
\begin{equation*} E_{m+1,n+1}=(-A^4-A^2)E_{m,n}\end{equation*}
It follows immediately that:
\begin{equation*}E_{m,n}=(-A^4-A^2)^n E_{m+n,0}\end{equation*}
Notice that $E_{m+n,0}$ is the projective line, or the projective line with 
an arrow on it.

We define a refinement of $<>_r$ with the help of the preceding formula:

\begin{definition}
Let $D$ be a diagram of a link in $\KB^1(\RR)$. Then:\\
$<D>_r=\Sum_i q_i E_{m_i,n_i}$ for some $q_i\in R$.

The {\it refined projective Kauffman bracket} of $D$, denoted 
$<>_{rp}$, is given by the formula:
\begin{equation*}<D>_{rp}=\Sum_i q_i (-A^4-A^2)^{n_i} E_{m_i+n_i,0}
\end{equation*}
\end{definition}

We may extend now Lemma \ref{lemma_xinout} to the case in which the vertical 
strand of the four diagrams of Figure~\ref{xinout} is in a projective 
line, if we replace $<>_r$ with $<>_{rp}$. We still assume that the 
diagrams have no crossings.

\begin{lemma}\label{lemma_xinout2}
The refined projective Kauffman bracket satisfies: 

(1) $<D_u>_{rp}=-A^{-2}<D_r>_{rp}-A^2<D_d>_{rp}$

(2) $<D_u>_{rp}=-A^{-4}<D_l>_{rp}-A^{-2}<D_d>_{rp}$
\begin{proof}[\rm{Proof}]
If the vertical strand belongs to an oval, this follows from Lemma \ref{lemma_xinout}
and the fact that $<>_{rp}$ is a refinement of $<>_r$.

Assume now, that the vertical strand is in a projective line. We have:\\
$<D_u>=<D_d>$ (the parity on the projective line does not change)\\
$<D_l>=<D_r>$ (in both cases $x$ is close to the projective line, so it
is not nested in another oval and it has the same contribution in $<D_l>$ 
and $<D_r>$).

Now, $<D_r>_r=E_{m,n}$, for some $m$ and $n$ and $<D_u>_r=E_{m-1,n-1}$. 
So:
\begin{equation*}<D_r>_{rp}=(-A^4-A^2)^n E_{m+n,0}\end{equation*}
\begin{equation*}<D_u>_{rp}=(-A^4-A^2)^{n-1}E_{m+n,0}\end{equation*}
Thus:
\begin{equation*}<D_r>_{rp}=(-A^4-A^2)<D_u>_{rp}=-A^4<D_d>_{rp}-A^2<D_u>_{rp}
\end{equation*}
Rearranging the terms gives (1).

Also:
\begin{equation*}<D_l>_{rp}=-A^4<D_u>_{rp}-A^2<D_d>_{rp}\end{equation*}
Rearranging the terms gives (2).
\end{proof}
\end{lemma}

We can now prove:

\begin{proposition}
The refined projective Kauffman bracket, $<>_{rp}$, is invariant under all
regular Reidemeister moves.

Thus, $\KB^1(\RR)$ is isomorphic with $R\oplus R$ and it is generated by
the projective line and the projective line with one arrow on it.
\begin{proof}[\rm{Proof}]
All that remains is the invariance of $<>_{rp}$ under $\Om_5$ moves. It is 
sufficient to prove it for diagrams with a single crossing, the one in 
the move. From Lemma~\ref{lemma_standard5} and Lemma~\ref{typesI_II} it 
follows that one has only to check the invariance of $<>_{rp}$ for a move 
that is of type IV in Figure~\ref{standard5}.
One may assume that there are only $x$-s inside and outside the loop containing the only
crossing, as all ovals are expressed with $x$-s already with $<>_r$.
Using Lemma~\ref{lemma_xinout2} it is possible to push all $x$-s inside 
the loop of the only crossing into the loop and all $x$-s outside the loop 
into the loop as well. After applying some $\Om_8$ moves to push all 
arrows on one part of the projective line, we obtain the move pictured in 
Figure~\ref{omega5}.

\begin{figure}[h]
\scalebox{0.8}{\includegraphics{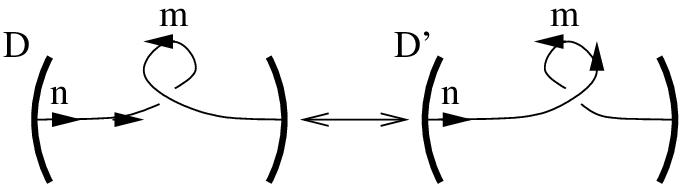}}
\caption{}
\label{omega5}
\end{figure}

Then, using Lemma~\ref{lemma_xinout2} again, it is possible to reduce $m$ 
to $0$ or $1$ by pushing some $x$-s out of the loop and into the left part 
of the projective line in Figure~\ref{omega5}.

For $m=0$: 
\begin{eqnarray*}
  <D'>_{rp} & = & A E_{n-1,0}+A^{-1}<E_{n,1}>_{rp}=A 
E_{n+1,0}+A^{-1}(-A^4-A^2)E_{n+1,0}= \\
  & = & -A^3 E_{n+1,0}=<D>_{rp}
  \end{eqnarray*}
the last equality from the framing relation.\\

For $m=1$: 
\begin{eqnarray*}
<D>_{rp}&=&A<E_{n+1,1}>_{rp}+A^{-1}E_{n,0}=A(-A^4-A^2)E_{n+2,0}+A^{-1}E_{n,0}=
\\ & = & (-A^5-A^3+A^{-1})E_{n,0}\\
<D'>_{rp} & = &A E_{n-2,0}+A^{-1}(-A^{-2}<E_{n,2}>_{rp}+(A^4+1)E_{n,0})
\end{eqnarray*}
because a circle with two arrows on it is $P_2$ and we have (see 
Definition \ref{P_n}):
\begin{equation*}<P_2>_r=-A^{-2}x<P_1>_r-A^2<P_0>_r=-A^{-2}x^2+A^4+1
\end{equation*}
So:
\begin{eqnarray*}
<D'>_{rp}&=&A E_{n,0}+A^{-1}(-A^{-2}(-A^4-A^2)^2E_{n,0}+(A^4+1)E_{n,0})=\\
& = & (A-A^{-3}(A^8+2A^6+A^4)+A^3+A^{-1})E_{n,0}=\\
& = & (-A^5-A^3+A^{-1})E_{n,0}=<D>_{rp}
\end{eqnarray*}
\end{proof}
\end{proposition}

\subsection{The submodule $\KB^0(\RR)$}

\begin{proposition}\label{rough_thm}
$\KB^0(\RR)$ is the quotient of $R[x]$ modulo relations:
\begin{equation*}(R_n):\; A P_{n-2}+A^{-1}P_{-n+2}=A P_{-n}+A^{-1}P_{n}, 
\;n\ge 2,\;n\in\N
\end{equation*}
where $P_n$ are polynomials in $R[x]$ from Definition \ref{P_n}.

\begin{proof}[\rm{Proof}]
For any diagram $D$ of a link in $\KB^0(\RR)$, $<D>_r\in R[x]$. 
From Lemma~\ref{lemma_standard5} and Lemma~\ref{typesI_II},
$<>_r$ is invariant under all regular Reidemeister moves except the 
$\Om_5$ moves of type III in Figure~\ref{standard5}.
Thus, $\KB^0(\RR)$ is the quotient of $R[x]$ modulo relations coming from 
$\Om_5$ moves of type III.

Using $\Omega_8$ moves and Lemma \ref{lemma_xinout} all arrows can be 
pushed on two strands, as in Figure~\ref{tors_5}.

\begin{figure}[h]
\scalebox{0.6}{\includegraphics{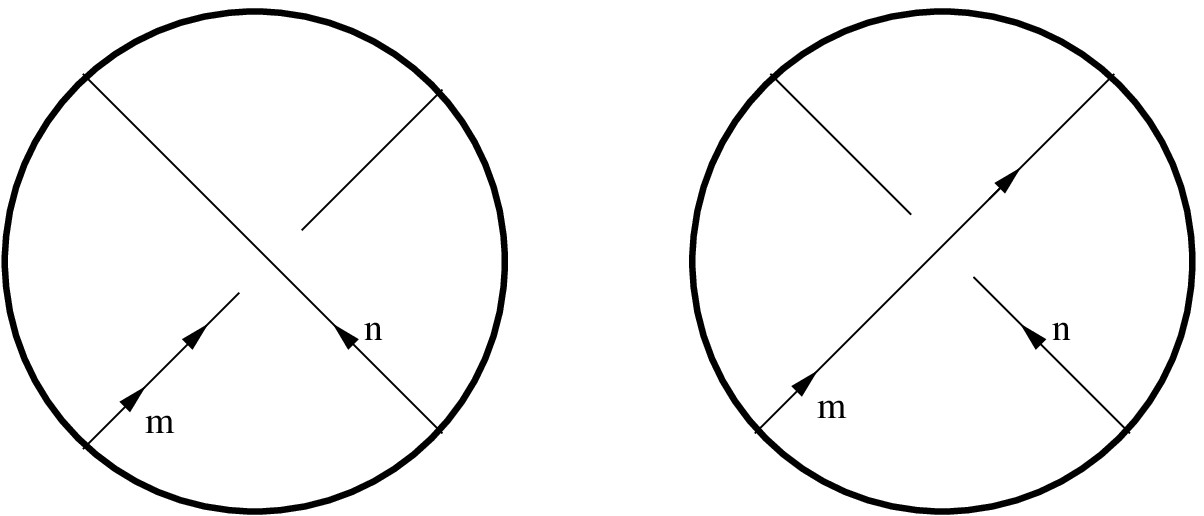}}
\caption{}
\label{tors_5}
\end{figure}

Applying Kauffman relation $(K1)$, followed by $\Om_6$ and $\Om_8$ moves,
now gives relations:
\begin{equation*}A P_{n-m-1}+A^{-1}P_{m+1-n}=A P_{m-n-1}+A^{-1}P_{n+1-m}, 
\; m,n\in\Z
\end{equation*}
Clearly, one can take $m=1$ without omitting any relation, which gives 
$(R_n)$, $n\in\Z$. Finally, noticing that $(R_n)$ is the same as $(R_{2-n})$ 
it is sufficient to take $(R_n)$ with $n\ge 2$ ($(R_1)$ is a trivial 
relation so it can be disregarded).
\end{proof}
\end{proposition}

In the rest of this subsection a simplification of $(R_n)$-s, $n\ge 2$, is 
presented.

A change of variable is useful for $P_n$-s: let $t=-A^{-3}x$ (as a diagram $t$ is $x$
with a negative kink added).

\begin{lemma}\label{sym}
For all $n\in\Z$, $P_{-n}(t,A)=P_n(t,A^{-1})$
\begin{proof}[\rm{Proof}]
We have:
\begin{equation*}P_0=-A^2-A^{-2},\; P_1=x=-A^3 t,\; P_{-1}=A^{-6}x=-A^{-3} 
t\end{equation*}
Thus:
\begin{equation*}P_0(t,A)=P_0(t,A^{-1}),\; P_{-1}(t,A)=P_1(t,A^{-1})
\end{equation*}
The recurrence relation for $P_n$ in $t$ is:
\begin{equation*}P_n=-A^{-2}x P_{n-1}-A^2 P_{n-2}=A t P_{n-1}-A^2 P_{n-2}
\end{equation*}
From the last relation one gets:
\begin{equation*}P_{n-2}=A^{-1} t P_{n-1}-A^{-2}P_n
\end{equation*}
Replacing $n$ with $2-n$: 
\begin{equation*}P_{-n}=A^{-1} t P_{-n+1}-A^{-2}P_{-n+2}
\end{equation*}

Suppose now, by induction, that $P_{-i}(t,A)=P_i(t,A^{-1})$ for $0\le 
i<n$.\\Then:
\begin{eqnarray*}
P_{-n}(t,A)&=&A^{-1} t P_{-n+1}(t,A)-A^{-2}P_{-n+2}(t,A)\\
 &=&A^{-1} t P_{n-1}(t,A^{-1})-A^{-2}P_{n-2}(t,A^{-1})\\
 &=&P_n(t,A^{-1})
\end{eqnarray*}
\end{proof}
\end{lemma}

Let $Q_n$, $n\in\N\cup\{-2,-1\}$ be defined by:
\begin{equation*}Q_{-2}=-1$, $Q_{-1}=0$ and $Q_n=t Q_{n-1}-Q_{n-2}
\end{equation*}
Thus $Q_0=1$ and $Q_1=t$.

\begin{lemma}\label{PQ}
For $n\ge 0$ we have:
\begin{equation*}P_n=-A^{n+2}Q_n+A^{n-2}Q_{n-2},
\;P_{-n}=-A^{-n-2}Q_n+A^{-n+2}Q_{n-2}\end{equation*}
\begin{proof}[\rm{Proof}]
Let:
\begin{equation*}P'_0=-A^2,\; P'_1=-A^3 t,\;P'_n=AtP'_{n-1}-A^2 
P'_{n-2}
\end{equation*}
\begin{equation*}P''_0=-A^{-2},\; P''_1=0,\; P''_n=AtP''_{n-1}-A^2 
P''_{n-2}
\end{equation*}
Then $P_n=P'_n+P''_n$.

Now it is easy to check by induction that:
\begin{equation*}P'_n=-A^{n+2}Q_n,\; n\ge 0
\end{equation*}
\begin{equation*}P''_n=A^{n-2}Q_{n-2},\; n\ge 2
\end{equation*}
Thus for $n\ge 2$: 
\begin{equation*}P_n=-A^{n+2}Q_n+A^{n-2}Q_{n-2}
\end{equation*}

For $n=0$ and $n=1$, $-A^{n+2}Q_n+A^{n-2}Q_{n-2}$ become respectively:
\begin{equation*}-A^2 Q_0+A^{-2}Q_{-2}=-A^2-A^{-2}=P_0 \end{equation*}
\begin{equation*}-A^3 Q_1+A^{-1}Q_{-1}=-A^3 t=P_1 \end{equation*}

The second part of the assertion follows from Lemma \ref{sym} and the fact
that $Q_n$ are polynomials in $t$ only.
\end{proof}
\end{lemma}

Recall that $(R_n)$, $n\ge 2$, is:
\begin{equation*}A P_{n-2}+A^{-1}P_{-n+2}=A P_{-n}+A^{-1}P_{n}
\end{equation*}
Using the preceding Lemma, it can be rewritten as:
\begin{equation*}A(-A^n Q_{n-2}+A^{n-4}Q_{n-4})+
A^{-1}(-A^{-n}Q_{n-2}+A^{-n+4}Q_{n-4})=\end{equation*}
\begin{equation*}A (-A^{-n-2}Q_n+A^{-n+2}Q_{n-2})+
A^{-1}(-A^{n+2}Q_n+A^{n-2}Q_{n-2})
\end{equation*}
Which, after rearranging the terms, becomes for $n\ge 2$:
\begin{equation*}(R'_n):\; 
(A^{n+1}+A^{-n-1})(Q_n-Q_{n-2})=(A^{n-3}+A^{-n+3})(Q_{n-2}-Q_{n-4}),
\end{equation*}

Let $q_{n,i}=\Sum_{k=1}^i (A^{n-4k+1}+A^{-n+4k-1})$

\begin{proposition}
For $1\le i\le \frac{n}{2}$, $(R'_n)$ can be rewritten using $(R'_k)$, 
$k<n$, as:
\begin{equation*}(A^{n+1}+A^{-n-1})(Q_n-Q_{n-2i})=q_{n,i}(Q_{n-2i}-Q_{n-2i-2})
\end{equation*}
\begin{proof}[\rm{Proof}]
The proof is by induction on $i$.
For $i=1$ this is $(R'_n)$.
Now suppose that the formula is true for $i\le \frac{n}{2}-1$. We have:
\begin{equation*}(A^{n+1}+A^{-n-1})(Q_n-Q_{n-2i-2})=\end{equation*}
\begin{equation*}(A^{n+1}+A^{-n-1})(Q_n-Q_{n-2i}+Q_{n-2i}-Q_{n-2i-2})=
\end{equation*}
\begin{equation*}(q_{n,i}+A^{n+1}+A^{-n-1})(Q_{n-2i}-Q_{n-2i-2})
\end{equation*}
the last equality by induction hypothesis.
\begin{equation*} 
(q_{n,i}+A^{n+1}+A^{-n-1})(Q_{n-2i}-Q_{n-2i-2})=(A^{n+1}+A^{n-3}+..
\end{equation*}
\begin{equation*}
+A^{n-4i+1}+A^{-n-1}+A^{-n+3}+..+A^{-n+4i-1})(Q_{n-2i}-Q_{n-2i-2})
\end{equation*}
Notice that:
\begin{equation*}(A^{n+1}+A^{-n+4i-1})=(A^{n-2i+1}+A^{-n+2i-1})A^{2i}
\end{equation*}
\begin{equation*}(A^{n-3}+A^{-n+4i-5})=(A^{n-2i+1}+A^{-n+2i-1})A^{2i-4}
\end{equation*}
\begin{equation*}...\end{equation*}
\begin{equation*}(A^{n-4i+1}+A^{-n-1})=(A^{n-2i+1}+A^{-n+2i-1})A^{2i-4i}
\end{equation*}
Thus:
\begin{equation*}(q_{n,i}+A^{n+1}+A^{-n-1})(Q_{n-2i}-Q_{n-2i-2})=
\end{equation*}
\begin{equation*}(A^{n-2i+1}+A^{-n+2i-1})
(Q_{n-2i}-Q_{n-2i-2})(A^{2i}+A^{2i-4}+..+A^{2i-4i})=
\end{equation*}
\begin{equation*}
(A^{n-2i-3}+A^{-n+2i+3})(Q_{n-2i-2}-Q_{n-2i-4})(A^{2i}+A^{2i-4}+..+A^{2i-4i})
\end{equation*}
where the last equality comes from $(R'_{n-2i})$, which can be applied as 
$n-2i\ge 2$.

So: 
\begin{equation*}(A^{n+1}+A^{-n-1})(Q_n-Q_{n-2i-2})= \end{equation*}
\begin{equation*}(A^{n-2i-3}+A^{-n+2i+3})
(A^{2i}+A^{2i-4}+..+A^{2i-4i})(Q_{n-2i-2}-Q_{n-2i-4})= \end{equation*}
\begin{equation*}
(A^{n-3}+A^{n-7}+..+A^{n-4i-3}+A^{-n+3}+A^{-n+7}+.. \end{equation*}
\begin{equation*}+A^{-n+4i+3})(Q_{n-2i-2}-Q_{n-2i-4})
=q_{n,i+1}(Q_{n-2i-2}-Q_{n-2i-4}) \end{equation*}
\end{proof}
\end{proposition}

\begin{corollary}
$(R'_n)$ can be rewritten as:

$(A^{n+1}+A^{-n-1})(Q_n-1)=2(A+A^{-1})\Sum_{k=1}^{\frac{n}{2}} 
A^{n+2-4k}$, for $n$ even,

$(A^{n+1}+A^{-n-1})(Q_n-t)=2t\Sum_{k=1}^{\frac{n-1}{2}} A^{n+1-4k}$, for 
$n$ odd.

\begin{proof}[\rm{Proof}]
If $n$ is even, one can take $i=\frac{n}{2}$ in the preceding Proposition.
Then:
\begin{equation*}Q_{n-2i}=Q_{0}=1,\; Q_{0}-Q_{-2}=1-(-1)=2 \end{equation*}
\begin{equation*}
q_{n,\frac{n}{2}}=A^{n-3}+A^{n-7}+..+A^{-n+1}+A^{-n+3}+A^{-n+7}+..+A^{n-1}=
\end{equation*}
\begin{equation*}
(A+A^{-1})(A^{n-2}+A^{n-6}+..+A^{-n+2})=(A+A^{-1})\Sum_{k=1}^{\frac{n}{2}} 
A^{n+2-4k} \end{equation*}

If $n$ is odd, one can take $i=\frac{n-1}{2}$ in the preceding Proposition.
Then:
\begin{equation*}Q_{n-2i}=Q_{1}=t,\; Q_1-Q_{-1}=t-0=t \end{equation*}
\begin{equation*} 
q_{n,\frac{n-1}{2}}=A^{n-3}+A^{n-7}+..+A^{-n+3}+A^{-n+3}+A^{-n+7}+..+A^{n-3}=
\end{equation*}
\begin{equation*}2\Sum_{k=1}^{\frac{n-1}{2}} A^{n+1-4k} \end{equation*}
\end{proof}
\end{corollary}

We can summarize the results of the two preceding subsections in the main 
theorem (note that we only need $Q_n$ with $n\ge 2$ to state this theorem):

\begin{theorem}\label{main_theorem}
$\KB(\RR)=R\oplus R\oplus R[t]/S$, where $R=\Z[A,A^{-1}]$ and $S$ is the submodule
of $R[t]$ generated by: 

$(A^{n+1}+A^{-n-1})(Q_n-1)-2(A+A^{-1})\Sum_{k=1}^{\frac{n}{2}}
A^{n+2-4k}$, for $n\ge 2$ even, 

$(A^{n+1}+A^{-n-1})(Q_n-t)-2t\Sum_{k=1}^{\frac{n-1}{2}} A^{n+1-4k}$, for 
$n\ge 3$ odd, 

where $Q_{0}=1$, $Q_1=t$ and $Q_n=t Q_{n-1}-Q_{n-2}$.
\end{theorem}

\begin{corollary}\label{can_form}
An element in $\KB(\RR)$ can be written in a unique way in the 
following form (called canonical):
\begin{equation*}
( r_1,r_2,p_0+p_1t+\Sum_{k=1}^n p_{2k}(Q_{2k}-1)+\Sum_{k=1}^m 
p_{2k+1}(Q_{2k+1}-t))\end{equation*}
where $n,m\in\N$, $r_1,r_2,p_0,p_1\in R$
and, for $k\ge 2$, $p_k\in\Z[A]$, $deg(p_k)\le 2k+1$
\begin{proof}[\rm{Proof}]
Clearly, an element in $\KB(\RR)$ can be written in the canonical form, 
using relations $(R'_n)$. Now, if two elements $m_1$ and $m_2$ are 
in the canonical form and represent the same element in $\KB(\RR)$, then\\ 
$m_1-m_2=(r_1,r_2,p_0+p_1t+\Sum_{k=1}^n p_{2k}(Q_{2k}-1)+\Sum_{k=1}^m
p_{2k+1}(Q_{2k+1}-t))$\\
is in the canonical form and 
belongs to the submodule $S$. Therefore $r_1=r_2=p_0=p_1=0$. Now, as 
$m_1-m_2\in S$, $p_k$ must be of the form $(A^{k+1}+A^{-k-1})p'_k$, for 
some $p'_k\in R$. As $deg(p_k)\le 2k+1$, this is possible only if 
$p'_k=0$. Thus, all $p_k$-s are zero and $m_1=m_2$. 
\end{proof}
\end{corollary}

We show now that there are torsion elements in $\KB(\RR)$.

\begin{proposition}
Let $n$ be even or $n=1\mbox{ }(mod\mbox{ }4)$, $n\ge 2$.
Then, there exists a torsion element $m\in\KB(\RR)$, such that 
$deg_t(m)=n$, where $m$ is in the canonical form.
\begin{proof}[\rm{Proof}]
Suppose that $n$ is even. Then $A^{n+1}+A^{-n-1}=(A+A^{-1})r_n$, where
$r_n=(A^n-A^{n-2}+..-A^{-n+2}+A^{-n})$. Let: \\
$m=A^n(r_n(Q_n-1)-2\Sum_{k=1}^{\frac{n}{2}}A^{n+2-4k})$\\
Then $A^n r_n$ is in $\Z[A]$ of degree $2n$, so $m$ is in canonical form 
and, from Corollary~\ref{can_form}, $m$ is not $0$ in $\KB(\RR)$.
However:\\
$(A+A^{-1})m=A^n((A^{n+1}+A^{-n-1})(Q_n-1)-2(A+A^{-1})
\Sum_{k=1}^{\frac{n}{2}}A^{n+2-4k})=0$\\
so $m$ is a torsion element.\\

Suppose now that $n=1\mbox{ }(mod\mbox{ }4)$.
Then $A^{n+1}+A^{-n-1}=(A^2+A^{-2})r_n$, where
$r_n=A^{n-1}-A^{n-5}+..-A^{-n+5}+A^{-n+1}$. Also:\\
$\Sum_{k=1}^{\frac{n-1}{2}} 
A^{n+1-4k}=(A^{n-3}+A^{n-7}+..+A^{-n+3})=
(A^2+A^{-2})(A^{n-5}+A^{n-13}+..$\\
$+A^{-n+5})=(A^2+A^{-2})\Sum_{k=1}^{\frac{n-1}{4}} A^{n+3-8k}$\\
Let:\\
$m=A^{n-1}(r_n(Q_n-t)-2t\Sum_{k=1}^{\frac{n-1}{4}} A^{n+3-8k})$\\
Then $A^{n-1} r_n$ is in $\Z[A]$ of degree $2n-2$, so $m$ is in canonical 
form and, from Corollary~\ref{can_form}, $m$ is not $0$ in $\KB(\RR)$.
However:\\
$(A^2+A^{-2})m=A^{n-1}((A^{n+1}+A^{-n-1})(Q_n-t)-2t(A^2+A^{-2})
\Sum_{k=1}^{\frac{n-1}{4}} A^{n+3-8k})=0$
so $m$ is a torsion element.
\end{proof} 
\end{proposition}

Now, we show that there are not always torsion elements in all degrees 
in the canonical form. 

Notice first that $R$ is a UFD. Indeed, the invertible elements in $R$ are 
of the form $\pm A^n$, and two different decompositions of an element in 
$R$ into irreducible elements would give, after multiplication by some 
invertible element, two different decompositions into irreducible elements 
in $\Z[A]$, a UFD. Also, an irreducible polynomial in $\Z[A]$ is 
irreducible in $R$. 

\begin{proposition}\label{no3}
Let $m\in\KB(\RR)$ be written in the canonical form. Suppose that 
$deg_t(m)=3$. Then $m$ is not a torsion element.
\begin{proof}[\rm{Proof}]
The polynomial $A^8+1$ is irreducible in $\Z[A]$, which can be checked by
replacing $A$ with $A+1$ and applying the Eisenstein irreducibility
criterion. Therefore, $A^4+A^{-4}$ is irreducible in $R$.

Now, the odd and even powers of $t$ are independent in $\KB(\RR)$ 
(they correspond to different homology classes in $H_1(\RR;\Z_2)$), so we 
may assume that $m=p_3(Q_3-t)+p_1 t$, with $p_3\in\Z[A]$, $deg(p_3)\le 7$, 
$p_3\neq 0$ and $p_1\in R$.

Suppose that for some $r\in R$, $r m=0$. Then $rp_3=(A^4+A^{-4})r'$, for 
some $r'\in R$. Now, as $(A^4+A^{-4})$ is irreducible in $R$, it divides 
$r$ or $p_3$. But it cannot divide $p_3$ as $deg(p_3)\le 7$. Thus 
$r=(A^4+A^{-4})w$, for some $w\in R$. We have then:
\begin{equation*}
rm=(A^4+A^{-4})w(p_3(Q_3-t)+p_1 t)=w(2+(A^4+A^{-4})p_1)t=0
\end{equation*}
Thus, $w(2+(A^4+A^{-4})p_1)=0$ and, as $R$ is a domain, $w=0$ or 
$2+(A^4+A^{-4})p_1=0$. The second case is ruled out by the fact that the 
lowest and highest powers in $A$ of $(A^4+A^{-4})p_1$ are different if 
$p_1\neq 0$. Thus $w=0$, so $r=0$ and $m$ is not a torsion element.
\end{proof}
\end{proposition}

Finally, unlike the case of $\KB(S^1\times S^2)$ (see Theorem 
\ref{thm_s2s1}), we have:

\begin{proposition}\label{non_split}
The module $\KB(\RR)$ does not split as a sum of cyclic $R$-modules.
\begin{proof}[\rm{Proof}]
Suppose that $\KB(\RR)=\Bigoplus_i R\oplus\Bigoplus_j R/I_j$, for 
some ideals $I_j$ in $R$. We can write: 
\begin{equation*}t=t_1+..+t_n+t_{tor}\end{equation*} 
where $t_1$,..,$t_n$ are non zero elements in different free $R$ summands 
of $\KB(\RR)$ and $t_{tor}$ is in $\Bigoplus_j R/I_j$ (as $t$ is not a 
torsion element $n$ must be at least $1$).

As $(A^4+A^{-4})(Q_3-t)=2t$, $2t_i=(A^4+A^{-4})r_i$, for some $r_i\in R$, 
$1\le i\le n$.
As $(A^4+A^{-4})$ is irreducible in $R$ (see the proof of Proposition 
\ref{no3}) and does not divide $2$, it divides $t_i$, so 
$t_i=(A^4+A^{-4})s_i$, for some $s_i\in R$.
Thus, for $s=s_1+..+s_n$:
\begin{equation*}t=(A^4+A^{-4})s+t_{tor}\end{equation*}
From Corollary~\ref{can_form}, it follows that for $m\in\KB(\RR)$, $m\neq 
0$ implies $2m\neq 0$. Thus, if for some $r\in R$, $2rm=0$, then $rm=0$.
Let $r\in R$, $r\neq 0$, be such that $rt_{tor}=0$. From the preceding 
remark, we may assume that $r\;(mod\;2)\neq 0$. Then:
\begin{equation*}rt=(A^4+A^{-4})rs\end{equation*}
If $s$, as a polynomial in $t$, contains some even powers of $t$, they 
must be killed by $(A^4+A^{-4})r$. Therefore, without changing the last 
equation, we may assume that $s$ contains only odd powers of $t$. So $s$ 
has the form:
\begin{equation*}s=p_1 t+p_3(Q_3-t)+..+p_k(Q_k-t),\;p_i\in 
R\end{equation*}
and:
\begin{equation*}rt=(A^4+A^{-4})(rp_1 t+rp_3(Q_3-t)+..+rp_k(Q_k-t))
\end{equation*}

For $i\ge 3$, all $(A^4+A^{-4})rp_i(Q_i-t)$ must be in $R t$, so 
$(A^4+A^{-4})rp_i(Q_i-t)-r'_i t\in S$, for some $r'_i\in R$ ($S$ is the 
submodule in Theorem \ref{main_theorem}). 
This implies that all $r'_i$ are divisible by $2$. So, for some $q\in R$:
\begin{eqnarray*}
rt &=&(A^4+A^{-4})rp_1 t+2q t\\
 &=&(A^4+A^{-4})rp_1 t\;(mod\;2)
\end{eqnarray*}
This is impossible from considerations of the highest and lowest power of 
$A$ in $r\;(mod\;2)$ and $(A^4+A^{-4})rp_1\;(mod\;2)$ and the fact that 
$r\;(mod\;2)\neq 0$.
\end{proof}
\end{proposition}

\section{The module $\KB(S^1\times S^2$)}

In \cite{HP} the skein module $\KB(S^1\times S^2)$ has 
been computed. We present in this section an alternative computation 
of this skein module, similar to the computation of $\KB^0(\RR)$.

\begin{figure}[h]
\scalebox{0.6}{\includegraphics{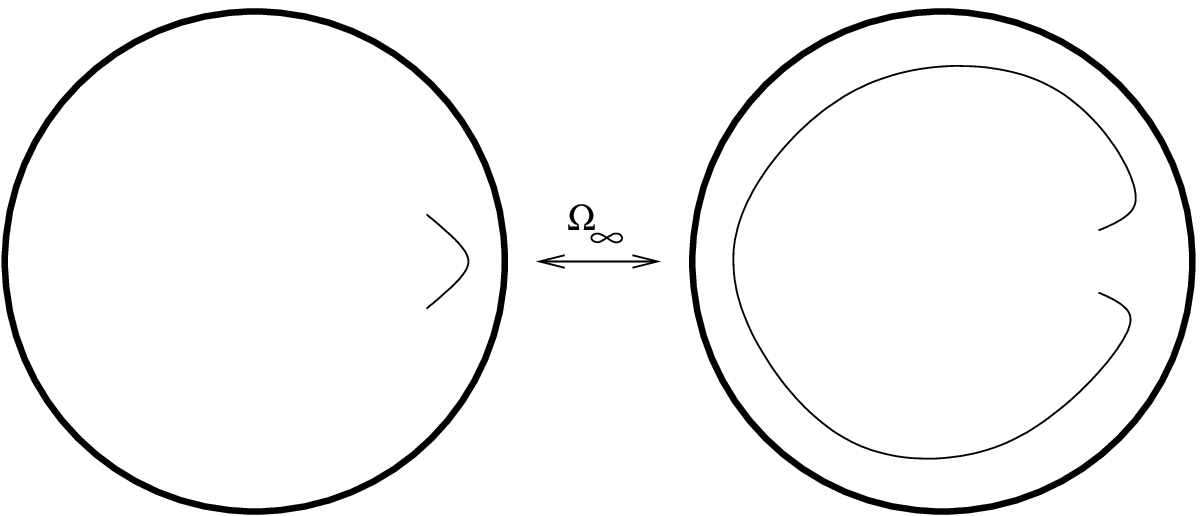}}
\caption{}
\label{om_inf}
\end{figure}

$S^1\times S^2$ is $S^1\times B^2/\partial{B^2}$. Therefore diagrams of 
links in $S^1\times S^2$ are the same as diagrams of links in $B^2\times 
S^1$. The additional generic singularity occurs when an arc contains the 
{\it infinity point} (i.e. $\partial{B^2}$ in $B^2/\partial{B^2}$).
Resolving it gives an additional $\Om_\infty$ move, presented in 
Figure~\ref{om_inf}.

One more time, it is possible to define a refined Kauffman bracket, 
$<>_r$, on diagrams of links in $S^1\times S^2$. Here it is exactly the 
$<>_r$ that was used to prove Proposition \ref{DS1}. In particular it is 
invariant under all regular Reidemeister moves except $\Om_\infty$ and it 
satisfies Lemma \ref{lemma_xinout}. 

\begin{proposition}\label{rough_thm2}
$\KB(S^1\times S^2)$ is the quotient of $R[x]$ modulo relations:
\begin{equation*}(S_n):\; P_{n}=P_{-n},\; n\ge 1,\; n\in\N
\end{equation*}
\begin{proof}[\rm{Proof}]
From the remarks before the Proposition it follows that $\KB(S^1\times 
S^2)$ is the quotient of $R[x]$ modulo all relations coming from  
$\Om_\infty$ moves. As before, we may consider such moves only on 
diagrams without crossings. Consider such a move, presented on the left of 
Figure~\ref{om_inf2}.

\begin{figure}[h]
\scalebox{0.6}{\includegraphics{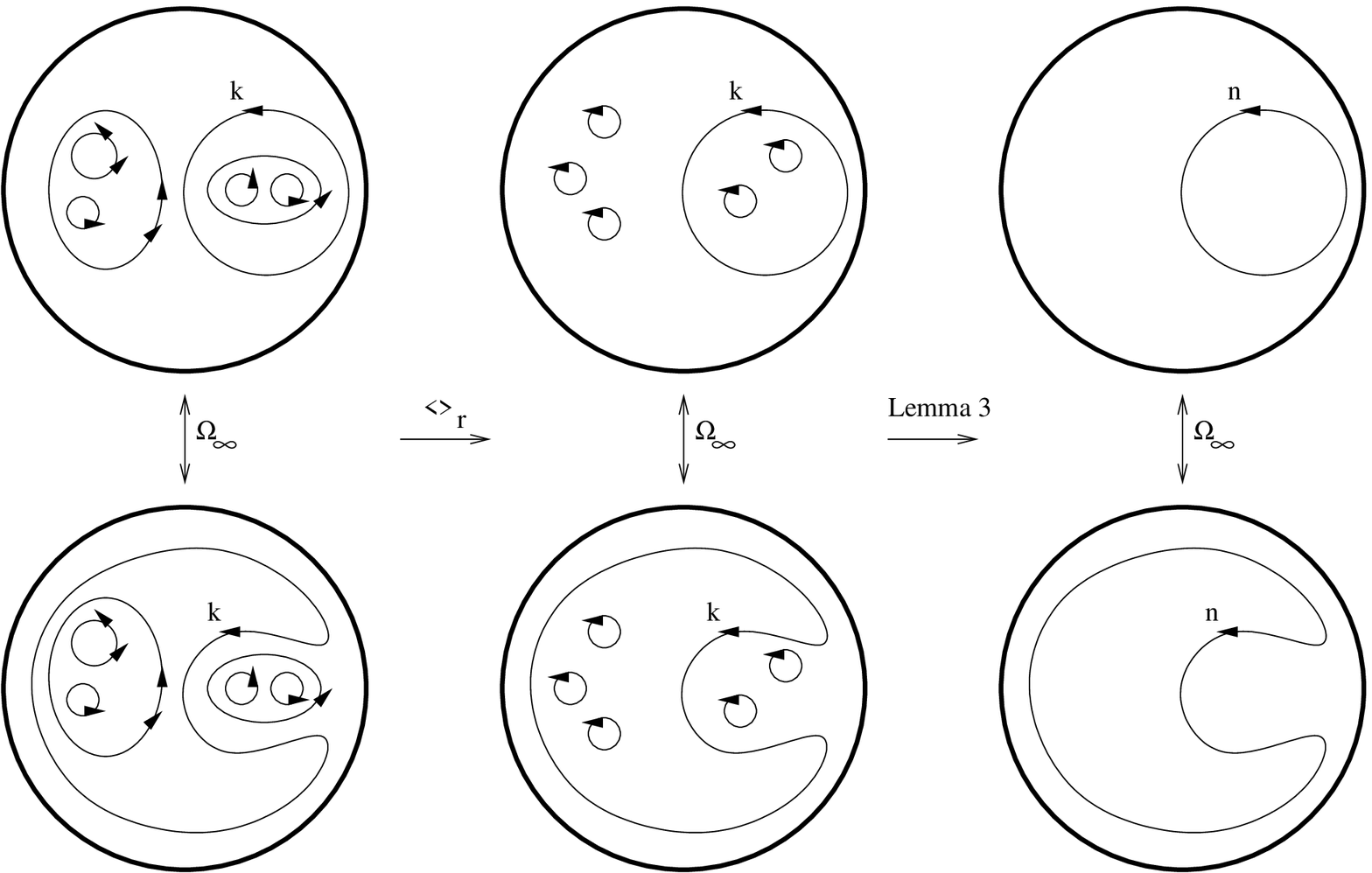}}
\caption{}
\label{om_inf2}
\end{figure}

We may apply $<>_r$ inside and outside the oval that is pushed through 
infinity, to reduce the situation to the middle of Figure~\ref{om_inf2}. 
Then, using Lemma \ref{lemma_xinout} all $x$-s are pushed inside this 
oval. Thus, the situation is simplified to the case presented on the 
right of Figure~\ref{om_inf2} in which a $P_n$ is pushed through 
infinity, becoming $P_{-n}$, $n\in\Z$. Thus, the relations needed to get
$\KB(S^1\times S^2)$ from $\KB(S^1\times B^2)$ are $P_n=P_{-n}$, $n\in 
\Z$. Obviously, it suffices to take $n\ge 1$.
\end{proof}
\end{proposition}

From Lemma \ref{PQ} the relations $(S_n)$ can be expressed as:
\begin{equation*}-A^{n+2}Q_n+A^{n-2}Q_{n-2}=-A^{-n-2}Q_n+A^{-n+2}Q_{n-2}
\end{equation*}
which can be rearranged becoming, for $n\ge 1$, the relations:
\begin{equation*}(S'_n):\;(-A^{n+2}+A^{-n-2})Q_n=(-A^{n-2}+A^{-n+2})Q_{n-2}
\end{equation*}

Let $\hat{q}_{n,i}=\Sum_{k=1}^i (-A^{n-4k+2}+A^{-n+4k-2})$

Let $Q_{n,i}=\Sum_{k=0}^i Q_{n-2k}$

\begin{proposition}
For $1\le i\le \frac{n+1}{2}$, $(S'_n)$ can be rewritten using $(S'_k)$, 
$k<n$, as:
\begin{equation*}(-A^{n+2}+A^{-n-2})Q_{n,i-1}=\hat{q}_{n,i}Q_{n-2i}
\end{equation*}

\begin{proof}[\rm{Proof}]
The proof is by induction on $i$.
For $i=1$ this is $(S'_n)$.
Now suppose that the formula is true for $i\le \frac{n-1}{2}$. We have:
\begin{equation*}
(-A^{n+2}+A^{-n-2})Q_{n,i}=(-A^{n+2}+A^{-n-2})(Q_{n,i-1}+Q_{n-2i})=
\end{equation*}
\begin{equation*}(\hat{q}_{n,i}-A^{n+2}+A^{-n-2})Q_{n-2i}\end{equation*}
the last equality by induction hypothesis.
\begin{equation*}(\hat{q}_{n,i}-A^{n+2}+A^{-n-2})Q_{n-2i}=\end{equation*}
\begin{equation*}(-A^{n+2}-A^{n-2}-A^{n-6}-..-A^{n-4i+2}+
A^{-n-2}+A^{-n+2}+..+A^{-n+4i-2})Q_{n-2i} \end{equation*}
Notice that:
\begin{equation*}
(-A^{n+2}+A^{-n+4i-2})=(-A^{n-2i+2}+A^{-n+2i-2})A^{2i}\end{equation*}
\begin{equation*}
(-A^{n-2}+A^{-n+4i-6})=(-A^{n-2i+2}+A^{-n+2i-2})A^{2i-4}\end{equation*}
\begin{equation*}...\end{equation*}
\begin{equation*}(-A^{n-4i+2}+A^{-n-2})=(-A^{n-2i+2}+A^{-n+2i-2})A^{2i-4i}
\end{equation*}
Thus:
\begin{equation*}(\hat{q}_{n,i}-A^{n+2}+A^{-n-2})Q_{n-2i}=\end{equation*}
\begin{equation*}
(-A^{n-2i+2}+A^{-n+2i-2})Q_{n-2i}(A^{2i}+A^{2i-4}+..+A^{2i-4i})=
\end{equation*}
\begin{equation*}
(-A^{n-2i-2}+A^{-n+2i+2})Q_{n-2i-2}(A^{2i}+A^{2i-4}+..+A^{2i-4i})
\end{equation*}
where the last equality comes from $(S'_{n-2i})$, which can be applied as 
$n-2i\ge 1$.

So: 
\begin{equation*}(-A^{n+2}+A^{-n-2})Q_{n,i}=\end{equation*}
\begin{equation*}(-A^{n-2i-2}+A^{-n+2i+2})
(A^{2i}+A^{2i-4}+..+A^{2i-4i})Q_{n-2i-2}=\end{equation*}
\begin{equation*}(-A^{n-2}-A^{n-6}-..-A^{n-4i-2}+
A^{-n+2}+A^{-n+6}+..+A^{-n+4i+2})Q_{n-2i-2}=\end{equation*}
\begin{equation*}\hat{q}_{n,i+1}Q_{n-2i-2}\end{equation*}
\end{proof}
\end{proposition}

Let $\hat{Q}_n=Q_{n,\frac{n}{2}}=Q_n+Q_{n-2}+..+Q_0$, for $n$ even\\ 
$\hat{Q}_n=Q_{n,\frac{n-1}{2}}=Q_n+Q_{n-2}+..+Q_1$, for $n$ odd.

\begin{theorem}\label{thm_s2s1}
$\KB(S^1\times S^2)=R\oplus\Bigoplus_{n=1}^\infty R/\{1-A^{2n+4}\}$
\begin{proof}[\rm{Proof}]
If $n$ is even, one can take $i=\frac{n}{2}$ in the preceding Proposition.
Then: 
\begin{equation*}\hat{q}_{n,\frac{n}{2}}=-A^{n-2}-A^{n-6}-..-A^{-n+2}
+A^{-n+2}+A^{-n+6}+..+A^{n-2}=0
\end{equation*}
Thus:
\begin{equation*}(-A^{n+2}+A^{-n-2})\hat{Q}_n=0
\end{equation*}

If $n$ is odd, one can take $i=\frac{n+1}{2}$ in the preceding Proposition.
Then:
\begin{equation*}
(-A^{n+2}+A^{-n-2})\hat{Q}_n=\hat{q}_{n,\frac{n+1}{2}}Q_{-1}=0
\end{equation*}
because $Q_{-1}=0$.

Multiplying these relations by invertible elements $A^{n+2}$, one gets:
\begin{equation*}(1-A^{2n+4})\hat{Q}_n=0
\end{equation*}

It suffices now to show that $\{\emptyset\}\cup\{\hat{Q}_n,n\ge 1\}$ is a 
basis of $\KB(S^1\times B^2)$. As $\{\emptyset\}\cup\{x^n,n\ge 1\}$ is 
such a basis, so $\{\emptyset\}\cup\{t^n,n\ge 1\}$ is also such a basis.
Now, $deg_t(\hat{Q}_n)=n$ and $\hat{Q}_n$ are monic polynomials in $t$, so 
they form, together with $\emptyset$, a basis of $\KB(S^1\times B^2)$.
\end{proof}
\end{theorem}

\section{The modules $\KB(L(p,1))$}

In \cite{HP2} the skein modules $\KB(L(p,q))$, $L(p,q)$ lens spaces, were 
computed. 
We present here a new computation of $\KB(L(p,1))$,
similar to the computation of $\KB(S^1\times S^2)$ of the previous 
section, though much simpler.

$L(p,1)$ is obtained from $S^1\times B^2$ by attaching a disk $B$ to 
it, where $\partial B$ is glued to a curve of type $(p,1)$ in 
$S^1\times\partial B^2$, then attaching a 3-ball. This last operation 
does not change the skein module, so we only need to consider the 
attaching of the disk.

\begin{figure}[h]
\scalebox{0.6}{\includegraphics{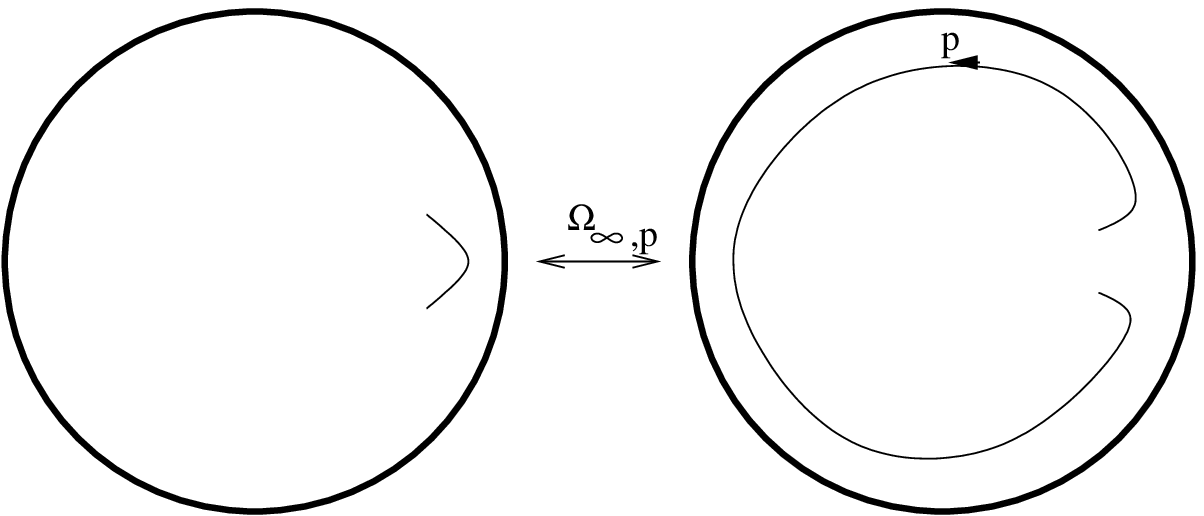}}
\caption{}
\label{om_inf_p}
\end{figure}

A link in $S^1\times B^2\cup B$ can be pushed along $B$ so 
that it lies in $S^1\times B^2$. Therefore, diagrams of links in $L(p,1)$ 
are the same as diagrams of links in $S^1\times B^2$. The extra 
Reidemeister move, $\Om_{\infty,p}$, comes from the sliding of an arc 
along $B$, so that it becomes a $(p,1)$ curve with a small 
segment removed. In a diagram, such a curve goes around the rest of the 
diagram and has $p$ arrows on it, see Figure~\ref{om_inf_p}.

It is possible to define a refined Kauffman bracket, $<>_r$, on diagrams 
of links in $L(p,1)$, exactly like in the previous section. So, again, it
is invariant under all regular Reidemeister moves except $\Om_{\infty,p}$ 
and it satisfies Lemma \ref{lemma_xinout}.

\begin{proposition}
$\KB(L(p,1))$ is the quotient of $R[x]$ modulo relations:
\begin{equation*}
(T_n):\; P_n=P_{p-n},\; n>\lfloor p/2 \rfloor,\; n\in \N
\end{equation*}
\begin{proof}[\rm{Proof}]
The proof follows the proof of Proposition \ref{rough_thm2}. Like before, 
$\KB(L(p,1))$ is the quotient of $R[x]$ modulo relations coming from 
$\Om_{\infty,p}$, and these relations can be reduced to the cases in which 
a $P_n$ is slided through the attached disk, becoming $P_{p-n}$, see 
Figure~\ref{om_inf_p2}. 

\begin{figure}[h]
\scalebox{0.6}{\includegraphics{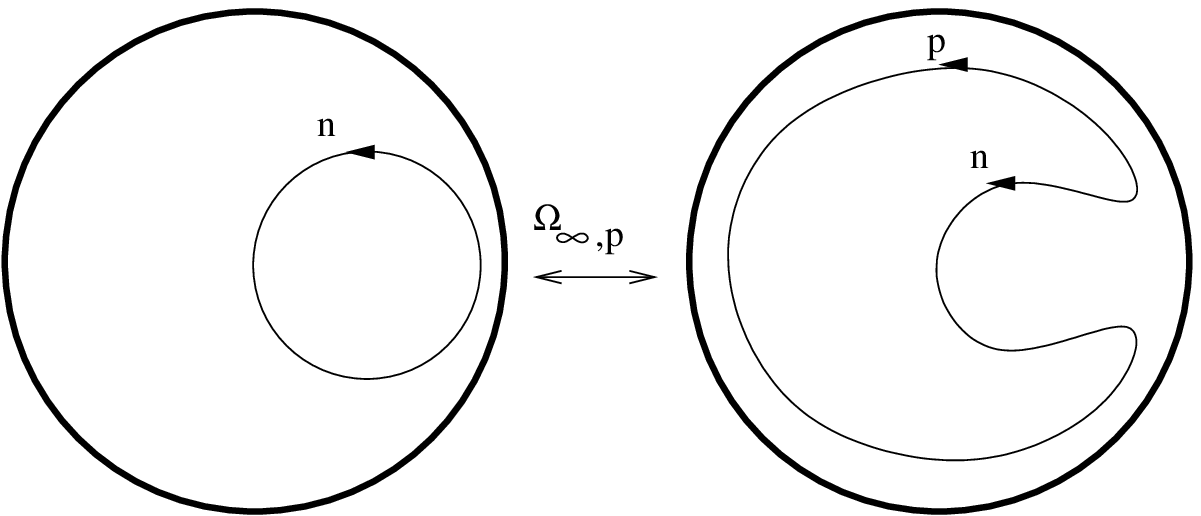}}
\caption{}
\label{om_inf_p2}
\end{figure}

The relation $(T_n)$ is obviously unchanged if one replaces $n$ with 
$p-n$, so it is sufficient to take $n>\lfloor p/2 \rfloor$ (for 
$p$ even $(T_{\frac{p}{2}})$ is a trivial relation so it can be 
disregarded).
\end{proof}
\end{proposition}

Note, that for $p=0$ this is just Proposition \ref{rough_thm2}, as 
$L(0,1)=S^1\times S^2$.

\begin{theorem}
Let $p\ge 1$. Then $\KB(L(p,1))=R[x]/(x^{\lfloor p/2 \rfloor+1})$.
\begin{proof}[\rm{Proof}]
From Definition \ref{P_n} it follows that, for $n>0$, 
$deg_x(P_n)=n$ and the leading coefficient of $P_n$ is an invertible 
element in $R$. From Lemma \ref{sym}, for $n<0$, 
$P_n(t,A)=P_{-n}(t,A^{-1})$, where $t=-A^{-3}x$. Thus 
$deg_x(P_n)=deg_t(P_n)=deg_t(P_{-n})=deg_x(P_{-n})=-n$. Also, as $t$ is 
equal to $x$ up to an invertible element in $R$, the leading coefficient 
of $P_n$ is an invertible element in $R$.
Thus, for all $n\in\Z$, $n\neq 0$, $P_n=a_n x^{|n|}+\hat{P}_n$, for some 
$a_n$ invertible in $R$ and some $\hat{P}_n$ 
satisfying $deg_x(\hat{P}_n)<|n|$.
Relations $(T_n)$, $n>\lfloor p/2\rfloor$, become:
\begin{equation*}a_n x^n+\hat{P}_n=a_{p-n}x^{|p-n|}+\hat{P}_{p-n}
\end{equation*}
Or:
\begin{equation*}x^n=a_n^{-1}(-\hat{P}_n+a_{p-n}x^{|p-n|}+\hat{P}_{p-n})
\end{equation*}
As, $p\ge 1$, $|p-n|<n$, so the degree in $x$ of the right hand side of 
the preceding relation is smaller than $n$.

Thus, these relations allow to express all $x^n$, 
$n>\lfloor p/2\rfloor$, with some lower degree polynomials. 
Therefore, as there are no other relations, $\KB(L(p,1))$ is free, 
generated by $\emptyset, x, ..., x^{\lfloor p/2\rfloor}$.
\end{proof}
\end{theorem}

The author was supported by the KBN grant KBN 0524/H03/2006/31.

\vspace{\baselineskip}
\begin{flushright}
\begin{tabular}{l}
Institute of Mathematics\\
University of Gdansk\\
ul. Wita Stwosza 57\\
80-952 Gdansk-Oliwa\\
Poland\\
e-mail: mmroczko@math.univ.gda.pl\\
\end{tabular}
\end{flushright}

\end{document}